
\def\ignore#1{}
 

\newcount\sectnum
\newcount\subsectnum
\newcount\eqnumber

\global\eqnumber=1\sectnum=0


\def\lab{(\the\sectnum.\the\eqnumber)}



\def\show#1{#1}



\def\smskip{\vskip 5 pt}
\def\medskip{\vskip 10 pt}
\def\bigskip{\vskip 15 pt}
\def\pn{\par\noindent}

\def\argmin{\mathop{\arg \min}}

\def\implies{\Rightarrow}

\def\frac#1#2{{#1\over #2}}

\def\ol#1{\overline{#1}}

\def\a{\alpha}

\def\m{\mu}

\def\re{\Re}

\def\tl{\tilde}

\def\old#1{}
\def\leaderfill{\leaders\hbox to 1em{\hss.\hss}\hfill}


\parindent=2pc
\baselineskip=15pt
\vsize=8.7 true in
\voffset=0.125 true in
\parskip=3pt


\def\minprob#1#2#3{$$\eqalign{&\hbox{minimize\ \ }#1\cr &\hbox{subject to\ \
}#2\cr}\ifnum 0=#3{}\else\eqno(#3)\fi$$}        
     
\def\maxprob#1#2#3{$$\eqalign{&\hbox{maximize\ \ }#1\cr &\hbox{subject to\ \
}#2\cr}\ifnum 0=#3{}\else\eqno(#3)\fi$$}        
     
\def\aligntwo#1#2#3#4#5{$$\eqalign{#1&#2\cr #3&#4\cr}
\ifnum 0=#5{}\else\eqno(#5)\fi$$}
\def\alignthree#1#2#3#4#5#6#7{$$\eqalign{#1&#2\cr #3&#4\cr #5&#6\cr}
\ifnum 0=#7{}\else\eqno(#7)\fi$$}


\def\eqnum{\eqno{\hbox{(\the\sectnum.\the\eqnumber)}\global\advance\eqnumber
by1}}

\def\eqnu{\eqno{\hbox{(\the\sectnum.\the\eqnumber)}\global\advance\eqnumber
by1}}

\newcount\examplnumber
\def\examplnum{\global\advance\examplnumber by1}

\newcount\figrnumber
\def\figrnum{\global\advance\figrnumber by1}

\newcount\propnumber
\def\propnum{\global\advance\propnumber by1}

\newcount\defnumber
\def\defnum{\global\advance\defnumber by1}

\newcount\lemmanumber
\def\lemmanum{\global\advance\lemmanumber by1}

\newcount\assumptionnumber
\def\assumptionnum{\global\advance\assumptionnumber by1}

\newcount\conditionnumber
\def\conditionnum{\global\advance\conditionnumber by1}

\def\exampl{\the\sectnum.\the\examplnumber}
\def\figr{\the\sectnum.\the\figrnumber}
\def\propn{\the\sectnum.\the\propnumber}
\def\defn{\the\sectnum.\the\defnumber}
\def\lemman{\the\sectnum.\the\lemmanumber}
\def\assumptionn{\the\sectnum.\the\assumptionnumber}
\def\condn{\the\sectnum.\the\conditionnumber}

\def\section#1{\goodbreak\vskip 3pc plus 6pt minus 3pt\leftskip=-2pc
   \global\advance\sectnum by 1\eqnumber=1
\global\examplnumber=1\figrnumber=1\propnumber=1\defnumber=1\lemmanumber=1\assumptionnumber=1 \conditionnumber =1%
   \line{\hfuzz=1pc{\hbox to 3pc{\bf 
   \vtop{\hfuzz=1pc\hsize=38pc\hyphenpenalty=10000\noindent\uppercase{\the\sectnum.\quad #1}}\hss}}
			\hfill}
			\leftskip=0pc\nobreak\tenf
			\vskip 1pc plus 4pt minus 2pt\noindent\ignorespaces}



\def\sect#1{\noindent\leftskip=-2pc\tenf
   \goodbreak\vskip 1pc plus 4pt minus 2pt
                \global\advance\subsectnum by 1\eqnumber=1
   \line{\hfuzz=1pc{\hbox to 3pc{\bf 
   \vtop{\hfuzz=1pc\hsize=38pc\hyphenpenalty=10000\noindent\uppercase{{\bf #1}}}\hss}}
                        \hfill}
   \leftskip=0pc\nobreak\tenf
                        \vskip 1pc plus 4pt minus 2pt\nobreak\noindent\ignorespaces}

\def\subsection#1{\noindent\leftskip=0pc\tenf
   \goodbreak\vskip 1pc plus 4pt minus 2pt
   \line{\hfuzz=1pc{\hbox to 3pc{\bf 
   \vtop{\hfuzz=1pc\hsize=38pc\hyphenpenalty=10000\noindent{\bf #1}}\hss}}
                        \hfill}
   \leftskip=0pc\nobreak\tenf
                        \vskip 1pc plus 4pt minus 2pt\nobreak\noindent\ignorespaces}
\def\subsubsection#1{\goodbreak\vskip 1pc plus 4pt minus 2pt
   \hfuzz=3pc\leftskip=0pc\noindent\tenit #1 \nobreak\tenf\vskip 6pt plus 1pt
                                minus 1pt\nobreak\ignorespaces\leftskip=0pc}
%

\def\beginexample#1{\noindent\goodbreak\vskip 6pt plus 1pt minus 1pt
\noindent
  \hbox {\bf Example #1\hss}
  \nobreak\vskip 4pt plus 1pt minus 1pt \nobreak\noindent\ninef
  \global\advance
                \leftskip by\parindent\pn}
\def\endexample{\vskip 12pt\tenf\par
  \global\advance\leftskip by -\parindent
  }

\def\beginexercise#1{\noindent\goodbreak\vskip 6pt plus 1pt minus 1pt \noindent\global\normalbaselineskip=12pt
  \hbox {\bf Exercise #1\hss}
  \nobreak\vskip 4pt plus 1pt minus 1pt 
  \nobreak\noindent\ninef\global\advance\leftskip
                        by\parindent\pn}
\def\endexercise{\vskip 12pt\tenf\par
  \global\advance\leftskip by -\parindent
  }

\def\beginsection#1{\noindent\goodbreak\vskip 6pt plus 1pt minus 1pt \noindent\global\normalbaselineskip=12pt
  \hbox {\it #1\hss}
  \vskip 0.1pt plus 1pt minus 1pt \nobreak\noindent\ninef\global\advance
                \leftskip by\parindent\noindent\pn}
\def\endsection{\vskip 12pt\tenf\par
  \global\advance\leftskip by -\parindent
}

%


\def\proposition#1{\smskip\pn{\bf Proposition #1}\quad}
\def\proof{\smskip\pn{\bf Proof:}\quad} 
\def\definition#1{\smskip\pn{\bf Definition #1}\quad} 
\def\assumption#1{\smskip\pn{\bf Assumption #1}\quad}

 \def\qed{\quad{\bf
Q.E.D.} \par\bigskip}
\def\ref{\smskip\pn}

\def\chapter#1#2{{\bf \centerline{\helbigbig
{#1}}}\bigskip\bigskip{\bf \centerline{\helbigbig
{#2}}}\bigskip\bigskip} 



\def\longpapertitle#1#2#3{{\bf \centerline{\helbigb
{#1}}}\bigskip{\bf \centerline{\helbigb
{#2}}}\bigskip\bigskip{\centerline{
by}}\bigskip{\bf \centerline{
{#3}}}\bigskip\bigskip} 


\def\nitem#1{\smskip\item{#1}}

\newcount\alphanum
\newcount\romnum

\def\alphaenumerate{\ifcase\alphanum \or (a)\or (b)\or (c)\or (d)\or (e)\or
(f)\or (g)\or (h)\or (i)\or (j)\or (k)\fi}
\def\romenumerate{\ifcase\romnum \or (i)\or (ii)\or (iii)\or (iv)\or (v)\or
(vi)\or (vii)\or (viii)\or (ix)\or (x)\or (xi)\fi}

\def\alist{\begingroup\vskip10pt\alphanum=1
\parskip=2pt\parindent=0pt \leftskip=3pc
\everypar{\llap{{\rm\alphaenumerate\hskip1em}}\advance\alphanum by1}}

\def\nolist{\begingroup\vskip10pt\alphanum=0
\parskip=2pt\parindent=0pt \leftskip=3pc
\everypar{\llap{\global\advance\alphanum by1(\the\alphanum)\hskip1em}}}

\def\romlist{\begingroup\vskip10pt\romnum=1
\parskip=2pt\parindent=0pt \leftskip=5pc
\everypar{\llap{{\rm\romenumerate\hskip1em}}\advance\romnum by1}}



\long\def\fig#1#2#3{\vbox{\vskip1pc\vskip#1
\prevdepth=12pt \baselineskip=12pt
\vskip1pc
\hbox to\hsize{\hfill\vtop{\hsize=25pc\noindent{\eightbf Figure #2\ }
{\eightpoint#3}}\hfill}}}

\long\def\widefig#1#2#3{\vbox{\vskip1pc\vskip#1
\prevdepth=12pt \baselineskip=12pt
\vskip1pc
\hbox to\hsize{\hfill\vtop{\hsize=28pc\noindent{\eightbf Figure #2\ }
{\eightpoint#3}}\hfill}}}

\long\def\table#1#2{\vbox{\vskip0.5pc
\prevdepth=12pt \baselineskip=12pt
\hbox to\hsize{\hfill\vtop{\hsize=25pc\noindent{\eightbf Table #1\ }
{\eightpoint#2}}\hfill}}}

 
\def\rightheadline#1{\headline{\tenrm\hfil #1}}


\long\def\leftfig#1#2{\vbox{\smskip\hsize=220pt
\vtop{{\noindent {\bf #1}}}
\smskip
\noindent
\vbox{{\noindent #2}}
}}

\long\def\rightfig#1#2#3{\vbox{\smskip\vskip#1
\prevdepth=12pt \baselineskip=12pt
\hsize=210pt
\smskip
\vbox{\noindent{\eightbold #2}
\hskip1em{\eightpoint#3}}
}}

\long\def\concept#1#2#3#4#5{\bigskip\hrule
\vbox{\hbox{\leftfig{#1}{#2} \hskip3em
\rightfig{#3}{#4}{#5}} \smskip}
\hrule\bigskip}


\long\def\bconcept#1#2#3#4#5#6#7{
\vbox{
\hbox to \hsize{\vtop{\par #1}}
\concept{#2}{#3}{#4}{#5}{#6}
\hbox to \hsize{\vtop{\par #7}}
\smskip}
}




\def\boxit#1{\vbox{\hrule\hbox{\vrule\kern3pt
                                \vbox{\kern3pt#1\kern3pt}\kern3pt\vrule}\hrule}}
\def\centerboxit#1{$$\vbox{\hrule\hbox{\vrule\kern3pt
                                \vbox{\kern3pt#1\kern3pt}\kern3pt\vrule}\hrule}$$}

\long\def\boxtext#1#2{$$\boxit{\vbox{\hsize #1\noindent\strut #2\strut}}$$}

%
%
%

\def\picture #1 by #2 (#3){
  \vbox to #2{
    \hrule width #1 height 0pt depth 0pt
    \vfill
    \special{picture #3} 
    }
  }

\def\scaledpicture #1 by #2 (#3 scaled #4){{
  \dimen0=#1 \dimen1=#2
  \divide\dimen0 by 1000 \multiply\dimen0 by #4
  \divide\dimen1 by 1000 \multiply\dimen1 by #4
  \picture \dimen0 by \dimen1 (#3 scaled #4)}
  }

%
%

\long\def\captfig#1#2#3#4#5{\vbox{\vskip1pc
\hbox to\hsize{\hfill{\picture #1 by #2 (#3)}\hfill}
\prevdepth=9pt \baselineskip=9pt
\vskip1pc
\hbox to\hsize{\hfill\vtop{\hsize=24pc\noindent{\eightbold Figure #4}
\hskip1em{\eightpoint#5}}\hfill}}}

%
%
%

\def\illustration #1 by #2 (#3){
  \vskip#2\hskip#1\special{illustration #3} 
    }

\def\scaledillustration #1 by #2 (#3 scaled #4){{
  \dimen0=#1 \dimen1=#2
  \divide\dimen0 by 1000 \multiply\dimen0 by #4
  \divide\dimen1 by 1000 \multiply\dimen1 by #4
  \illustration \dimen0 by \dimen1 (#3 scaled #4)}
  }


\newbox\graybox
\newdimen\xgrayspace
\newdimen\ygrayspace
%
%
%
%
%
%
%
%
%

\def\Textshade#1#2#3#4#5#6{%
    \xgrayspace=#4pt%
    \ygrayspace=#4pt%
    \def\grayshade{#3}%
    \def\linewidth{#5}%
    \def\theradius{#6}%
    \setbox\graybox=\hbox{\surroundboxa{#2}}%
    \hbox{%
    \hbox to 0pt{%
    \PScommands
}%
    \box\graybox}}%
%
%
\long%

\long%
\def\Parashade#1#2#3#4#5#6#7{%
    \xgrayspace=#4pt%
    \ygrayspace=#4pt%
    \def\grayshade{#3}%
    \def\linewidth{#5}%
    \def\theradius{#6}%
    \def\thevskip{#7pt}%
    \setbox\graybox=\hbox{\surroundboxb{#2}}%
    \vskip\thevskip%
    \hbox{%
    \hbox to 0pt{%
    \PScommands
}%
     \box\graybox}%
     \vskip\thevskip%
}%
%
%
%
\long\def\surroundboxa#1{\leavevmode\hbox{\vtop{%
\vbox{\kern\ygrayspace%
\hbox{\kern\xgrayspace#1%
      \kern\xgrayspace}}\kern\ygrayspace}}}
%
%
\long\def\surroundboxb#1{\leavevmode\hbox{\vtop{%
\vbox{\kern\ygrayspace%
\hbox{\kern\xgrayspace\vbox{\advance\hsize-2\xgrayspace#1}%
      \kern\xgrayspace}}\kern\ygrayspace}}}
%
%
%
\long\def\PScommands{%
\special{rawpostscript
/sharpbox{%
           newpath
           xmin ymin moveto
           xmin ymax lineto
           xmax ymax lineto
           xmax ymin lineto
           xmin ymin lineto
           closepath 
          }bind def
}%
\special{rawpostscript
/sharpboxnb{%
           newpath
           xmin ymin moveto
           xmin ymax lineto
           xmax ymax lineto
           xmax ymin lineto
          }bind def
}%
\special{rawpostscript
/sharpboxnt{%
           newpath
           xmin ymax moveto
           xmin ymin lineto
           xmax ymin lineto
           xmax ymax lineto
          }bind def
}%
\special{rawpostscript
/roundbox{%
           newpath
           xmin radius add ymin moveto
           xmax ymin xmax ymax radius arcto
           xmax ymax xmin ymax radius arcto
           xmin ymax xmin ymin radius arcto
           xmin ymin xmax ymin radius arcto 16 {pop} repeat
           closepath
          }bind def
}%
\special{rawpostscript
/sharpcorners{%
               sharpbox gsave grayshade setgray fill grestore 
               linewidth setlinewidth stroke
              }bind def
}%
\special{rawpostscript
/sharpcornersnt{%
               sharpboxnt gsave grayshade setgray fill grestore 
               linewidth setlinewidth stroke
              }bind def
}%
\special{rawpostscript
/sharpcornersnb{%
               sharpboxnb gsave grayshade setgray fill grestore 
               linewidth setlinewidth stroke
              }bind def
}%
\special{rawpostscript
/roundcorners{%
               roundbox gsave grayshade setgray fill grestore 
               linewidth setlinewidth stroke
              }bind def
}%
\special{rawpostscript
/plainbox{%
           sharpbox grayshade setgray fill 
          }bind def
}%
%
\special{rawpostscript
/roundnoframe{%
               roundbox grayshade setgray fill 
              }bind def
}%
\special{rawpostscript
/sharpnoframe{%
               sharpbox grayshade setgray fill 
              }bind def
}%
}%
%
%

\def\pshade#1{\Parashade{sharpcorners}{#1}{0.95}{10}{0.5}{10}{10}}


\def\boxit#1{\vbox{\hrule\hbox{\vrule\kern3pt
                                \vbox{\kern3pt#1\kern3pt}\kern3pt\vrule}\hrule}}

\def\boxitnb#1{\vbox{\hrule\hbox{\vrule\kern3pt
                                \vbox{\kern3pt#1\kern3pt}\kern3pt\vrule}}}

\def\boxitnt#1{\vbox{\hbox{\vrule\kern3pt
                                \vbox{\kern3pt#1\kern3pt}\kern3pt\vrule}\hrule}}

\long\def\boxtext#1#2{$$\boxit{\vbox{\hsize #1\noindent\strut #2\strut}}$$}



\def\texshopbox#1{\boxtext{462pt}{\vskip-1.5pc\pshade{\vskip-1.0pc#1\vskip-2.0pc}}}


%
%
%
%
%
%
%
%
\font\helbigbig=cmr10 scaled 2500%
\font\helbigb=cmbx10 scaled 1500%
\font\eightbold=cmbx8%

\def\tenf{\hel}%
\def\tenit{\heli}%
\def\ninef{\ninehel}%
\def\nineit{\nineheli}%
%
%


\font\tenrm=cmr10%
\font\teni=cmmi10%
\font\tensy=cmsy10%
\font\tenbf=cmbx10%
\font\tentt=cmtt10%
\font\tenit=cmti10%
\font\tensl=cmsl10%

\def\tenpoint{\def\rm{\fam0\tenrm}%
\textfont0=\tenrm%
\textfont1=\teni%
\textfont2=\tensy%
\textfont\itfam=\tenit%
\textfont\slfam=\tensl%
\textfont\ttfam=\tentt%
\textfont\bffam=\tenbf%
\scriptfont0=\sevenrm%
\scriptfont1=\seveni%
\scriptfont2=\sevensy%
\scriptscriptfont0=\sixrm%
\scriptscriptfont1=\sixi%
\scriptscriptfont2=\sixsy%
\def\it{\fam\itfam\tenit}%
\def\tt{\fam\ttfam\tentt}%
\def\sl{\fam\slfam\tensl}%
\scriptfont\bffam=\sevenbf%
\scriptscriptfont\bffam=\sixbf%
\def\bf{\fam\bffam\tenbf}%
\normalbaselineskip=18pt%
\normalbaselines\rm}%

\font\ninerm=cmr9%
\font\ninebf=cmbx9%
\font\nineit=cmti9%
\font\ninesy=cmsy9%
\font\ninei=cmmi9%
\font\ninett=cmtt9%
\font\ninesl=cmsl9%

\def\ninepoint{\def\rm{\fam0\ninerm}%
\textfont0=\ninerm%
\textfont1=\ninei%
\textfont2=\ninesy%
\textfont\itfam=\nineit%
\textfont\slfam=\ninesl%
\textfont\ttfam=\ninett%
\textfont\bffam=\ninebf%
\scriptfont0=\sixrm%
\scriptfont1=\sixi%
\scriptfont2=\sixsy%
\def\it{\fam\itfam\nineit}%
\def\tt{\fam\ttfam\ninett}%
\def\sl{\fam\slfam\ninesl}%
\scriptfont\bffam=\sixbf%
\scriptscriptfont\bffam=\fivebf%
\def\bf{\fam\bffam\ninebf}%
\normalbaselineskip=16pt%
\normalbaselines\rm}%

\font\eightrm=cmr8%
\font\eighti=cmmi8%
\font\eightsy=cmsy8%
\font\eightbf=cmbx8%
\font\eighttt=cmtt8%
\font\eightit=cmti8%
\font\eightsl=cmsl8%

\def\eightpoint{\def\rm{\fam0\eightrm}%
\textfont0=\eightrm%
\textfont1=\eighti%
\textfont2=\eightsy%
\textfont\itfam=\eightit%
\textfont\slfam=\eightsl%
\textfont\ttfam=\eighttt%
\textfont\bffam=\eightbf%
\scriptfont0=\sixrm%
\scriptfont1=\sixi%
\scriptfont2=\sixsy%
\scriptscriptfont0=\fiverm%
\scriptscriptfont1=\fivei%
\scriptscriptfont2=\fivesy%
\def\it{\fam\itfam\eightit}%
\def\tt{\fam\ttfam\eighttt}%
\def\sl{\fam\slfam\eightsl}%
\scriptscriptfont\bffam=\fivebf%
\def\bf{\fam\bffam\eightbf}%
\normalbaselineskip=14pt%
\normalbaselines\rm}%

\font\sevenrm=cmr7%
\font\seveni=cmmi7%
\font\sevensy=cmsy7%
\font\sevenbf=cmbx7%

\def\sevenpoint{%
   \def\rm{\sevenrm}\def\bf{\sevenbf}%
   \def\smc{\sevensmc}\baselineskip=12pt\rm}%

\font\sixrm=cmr6%
\font\sixi=cmmi6%
\font\sixsy=cmsy6%
\font\sixbf=cmbx6%

\fontdimen13\tensy=2.6pt%
\fontdimen14\tensy=2.6pt%
\fontdimen15\tensy=2.6pt%
\fontdimen16\tensy=1.2pt%
\fontdimen17\tensy=1.2pt%
\fontdimen18\tensy=1.2pt%

\def\tenf{\tenpoint}%
\def\ninef{\ninepoint}%
%



\def\section#1{\goodbreak\vskip 3pc plus 6pt minus 3pt\leftskip=-2pc
   \global\advance\sectnum by 1\eqnumber=1\subsectnum=0%
\global\examplnumber=1\figrnumber=1\propnumber=1\defnumber=1\lemmanumber=1\assumptionnumber=1 \conditionnumber =1%
   \line{\hfuzz=1pc{\hbox to 3pc{\bf 
   \vtop{\hfuzz=1pc\hsize=38pc\hyphenpenalty=10000\noindent\uppercase{\the\sectnum.\quad #1}}\hss}}
			\hfill}
			\leftskip=0pc\nobreak\tenf
			\vskip 1pc plus 4pt minus 2pt\noindent\ignorespaces}
\def\subsection#1{\noindent\leftskip=0pc\tenf
   \goodbreak\vskip 1pc plus 4pt minus 2pt
               \global\advance\subsectnum by 1
   \line{\hfuzz=1pc{\hbox to 3pc{\bf \the\sectnum.\the\subsectnum.
   \vtop{\hfuzz=1pc\hsize=38pc\hyphenpenalty=10000\noindent{\bf #1}}\hss}}
                        \hfill}
   \leftskip=0pc\nobreak\tenf
                        \vskip 1pc plus 4pt minus 2pt\nobreak\noindent\ignorespaces}



\def\texshopbox#1{\boxtext{462pt}{\vskip-1.5pc\pshade{\vskip-1.0pc#1\vskip-2.0pc}}}


\input miniltx

\ifx\pdfoutput\undefined
  \def\Gin@driver{dvips.def} 
\else
  \def\Gin@driver{pdftex.def} 
\fi

\input graphicx.sty
\resetatcatcode

\long\def\fig#1#2#3{\vbox{\vskip1pc\vskip#1
\prevdepth=12pt \baselineskip=12pt
\vskip1pc
\hbox to\hsize{\hfill\vtop{\hsize=30pc\noindent{\eightbf Figure #2\ }
{\eightpoint#3}}\hfill}}}

\def\show#1{}

\rightheadline{\botmark}

\pageno=1

\def\longpapertitle#1#2#3{{\bf \centerline{\helbigb
{#1}}}\medskip{\bf \centeline{\helbigb
{#2}}}\bigskip{\bf \centerline{
{#3}}}\bigskip}

\vskip-3pc

\def\xstar{X^{\raise0.04pt\hbox{\sevenpoint *}} }

\def\jstar{J^{\raise0.04pt\hbox{\sevenpoint *}} }
\def\qstar{Q^{\raise0.04pt\hbox{\sevenpoint *}} }

\rightheadline{\botmark}

\pageno=1

\rightheadline{\botmark}

\pageno=1

\rightheadline{\botmark}

\pn {\bf April 2020}
\bigskip \bigskip\bigskip

\bigskip

\def\longpapertitle#1#2#3{{\bf \centerline{\helbigb
{#1}}}\medskip{\bf \centerline{\helbigb
{#2}}}\bigskip{\bf \centerline{
{#3}}}\bigskip}

\vskip-3pc

\longpapertitle{Multiagent Value Iteration Algorithms in }{Dynamic Programming and Reinforcement Learning}{ {Dimitri Bertsekas\footnote{\dag}{\ninepoint McAfee Professor of Engineering, MIT, Cambridge, MA, and Fulton Professor of Computational Decision Making, ASU, Tempe, AZ.}}}



\centerline{\bf Abstract}

\smskip
\pn We consider infinite horizon dynamic programming problems, where the control at each stage consists of several distinct decisions, each one made by one of several agents. In an earlier work we introduced a policy iteration algorithm, where the policy improvement is done one-agent-at-a-time in a given order, with knowledge of the choices of the preceding agents in the order. As a result, the amount of computation for each policy improvement grows linearly with the number of agents, as opposed to exponentially for the standard all-agents-at-once method. For the case of a finite-state discounted problem, we  showed convergence to  an agent-by-agent optimal policy. In this paper, this result is extended to value iteration and optimistic versions of policy iteration, as well as to more general DP problems where the Bellman operator is a contraction mapping, such as stochastic shortest path problems with all policies being proper.

\vskip-1pc

\section{Multiagent Problem Formulation}

\pn We consider an abstract form of infinite horizon dynamic programming (DP) problem, which contains as special case finite-state discounted Markovian decision problems (MDP), as well as more general problems where the Bellman operator is a monotone weighted sup-norm contraction. The distinguishing feature of the problem is that the control $u$ consists of $m$ components $u_\ell$, $\ell=1,\ldots,m$, where $m>1$:
$$u=(u_1,\ldots,u_m).\xdef\controlstruct{\lab}\eqnum\show{oneo}$$
Conceptually, each component may be viewed as chosen by a distinct agent, with knowledge of the selections of the other agents. We consider value iteration (VI) algorithms that involve minimization component-by-component as opposed to minimization over all components at once. This is similar to what is done in coordinate descent methods for multivariable optimization, and can lead to dramatic gains in computational efficiency for large and even moderate values of $m$. We propose several methods and we show their convergence to an agent-by-agent optimal policy, a type of policy that is related to the notion of person-by-person optimality from the theory of teams. Our analysis extends and complements our earlier proposals of rollout and policy iteration (PI) algorithms [Ber20].

We assume that $u$ is chosen from a finite constraint set $U(x)$ when the system is at state $x$. In our earlier paper [Ber20], we have made a stronger assumption: we assumed that each control component $u_\ell$, $\ell=1,\ldots,m$, is separately constrained to lie in a given finite set $U_\ell (x)$. 
In this case $U(x)$ is the Cartesian product set
$$U(x)=U_1(x)\times\cdots\times U_m(x).\xdef\cartprod{\lab}\eqnum\show{oneo}$$
In this paper, we do not impose this assumption, except occasionally to discuss its implications. As a result our algorithms must ensure that the selection of a control component at a given state and stage does not preclude the feasibility of selection of the other control components at the same state and stage. This complicates our algorithms relative to the Cartesian product case \cartprod. We will discuss the mechanism for dealing with this issue in Section 2. For the remainder of this section, we will assume no special structure for the constraint set $U(x)$ other than finiteness.

\subsubsection{The $\a$-Discounted MDP Case}

\pn A major context for application of our algorithmic ideas is the standard infinite horizon discounted MDP with states $x=1,\ldots,n$. Here, at state $x$, a control $u$ is applied, and the system transitions to a next state $y$ with transition probabilities $p_{xy}(y)$ and cost $g(x,u,y)$. The control is chosen at state $x$ from a finite constraint set $U(x)$. The cost function of a stationary policy $\m$ that applies control $\m(x)\in U(x)$ at state $x$ is denoted by $J_\m(x)$, and the optimal cost [the minimum over $\m$ of $J_\m(x)$] is denoted by $\jstar(x)$. 

The standard VI algorithm starts from some initial guess $J^0$ and iterates as follows:\footnote{\dag}{\ninepoint Throughout the paper, we will be using componentwise equality and inequality notation, whereby for any pair of real-valued functions $J,J'$, we write $J= J'$ (or $J\le J'$) if $J(x)=J'(x)$ [or  $J(x)\le J'(x)$, respectively] for all $x\in X$.}

$$J^{k+1}=TJ^k,\qquad k=0,1,\ldots,$$
where $T$ is the Bellman operator, which maps a vector $J=\big(J(1),\ldots,J(n)\big)$ to the vector $$TJ=\big((TJ)(1),\ldots,(TJ)(n)\big)$$
 according to
$$(TJ)(x)=\min_{u\in U(x)}\sum_{y=1}^np_{xy}(u)\big(g(x,u,y)+\a J(y)\big),\qquad x=1,\ldots,n.\xdef\tmap{\lab}\eqnum\show{oneo}$$
Thus each VI involves a comparison of all the Q-factors
 $$Q(x,u)=\sum_{y=1}^np_{xy}(u)\big(g(x,u,y)+\a J(y)\big),\qquad x=1,\ldots,n,\ u\in U(x).\xdef\qfactor{\lab}\eqnum\show{oneo}$$

A related algorithm is {\it optimistic PI\/}, which involves simultaneous value and policy iterations, using the Bellman operator $T_\m$ defined for each policy $\m$ by
$$(T_\m J)(x)=\sum_{y=1}^np_{xy}\big(\m(x)\big)\Big(g\big(x,\m(x),y\big)+\a J(y)\Big),\qquad x=1,\ldots,n.\xdef\tmumap{\lab}\eqnum\show{oneo}$$
Given a pair $(\m^k,J^k)$, this algorithm generates $(\m^{k+1},J^{k+1})$ according to 
$$T_{\m^{k+1}}J^k=TJ^k,\qquad J^{k+1}=T_{\m^{k+1}}^{q}J^k,\qquad k=0,1,\ldots,\xdef\optpi{\lab}\eqnum\show{oneo}$$
where $q$ is a positive integer (which in some cases may depend on $k$), and $T_\m^q$ denotes the mapping obtained by $q$-fold application of the mapping $T_\m$. When $q=1$ we obtain the VI algorithm $J^{k+1}=TJ^k$ and when $q\to\infty$, we have $J^{k+1}=J_{\m^k}$ (in the limit), so the algorithm approaches the standard PI algorithm where $\m^{k+1}$ is obtained from $\m^k$ according to
$$T_{\m^{k+1}}J_{\m^k}=TJ_{\m^k}.\xdef\standardpi{\lab}\eqnum\show{oneo}$$
 
Unfortunately, iterating with the mapping $T$ is inconvenient for problems involving even a moderate number of agents, because the size of the control constraint set $U(x)$ typically grows exponentially with $m$. In particular, in the Cartesian product case \cartprod, if each constraint set $U_\ell(x)$ consists of at most $s$ elements, minimization over $U(x)$ involves a comparison of as many as $s^m$ Q-factors of the form \qfactor. This motivates us to consider versions of the preceding algorithms that involve a simpler form of minimization. For example, minimization  over the component constraint sets $U_\ell(x)$, one component at a time, which involves a comparison of $s$ Q-factors for each agent, for a total of $s\cdot m$ Q-factors. 

\subsubsection{The General Contractive DP Case}

\pn It is convenient and useful to develop our algorithm in a more general setting, which  involves an  operator-based framework from the author's abstract DP book [Ber18]. In particular, we consider {\it a finite set $X$ of states  and a finite set $U$ of controls\/}, and for each $x\in X$, a nonempty control constraint set $U(x)\subset U$.\footnote{\dag}{\ninepoint The abstract DP framework of [Ber18] does not require finiteness of the state and control spaces. We impose the finiteness assumption in order to obtain the most powerful algorithmic results possible. However, at several points in the paper, and particularly in Section 5, we speculate around the possibility of extending our algorithms and analysis to infinite state and control spaces.}
We denote by ${\cal M}$ the set of all functions $\m:X\mapsto U$ with $\m(x)\in U(x)$ for all $x\in X$, which we refer to as {\it policies\/}. We introduce
a mapping $H:X\times U\times {\cal R}(X)\mapsto\re$, where $\re$ denotes the real line and ${\cal R}(X)$ denotes the set of real-valued functions $J:X\mapsto\re$. 
For each policy $\m\in {\cal M}$, we consider the mapping $T_\m:{\cal R}(X)\mapsto {\cal R}(X)$ defined by
$$(T_\m J)(x)=H\big(x,\m(x),J\big),\qquad x\in X.$$
We also consider the mapping $T$ defined by
$$(TJ)(x)=\min_{u\in U(x)}H(x,u,J)=\min_{\m\in{\cal M}}(T_\m J)(x),\qquad x\in X.$$
Note that the $\a$-discounted MDP is obtained when $H$ is given by
$$H(x,u,J)=\sum_{y=1}^np_{xy}(u)\big(g(x,u,y)+\a J(y)\big),\qquad x=1,\ldots,n.\xdef\mdph{\lab}\eqnum\show{oneo}$$

The problem is to find a function $\jstar \in {\cal R}(X)$ such that
$$\jstar (x)=\min_{u\in U(x)}H(x,u,\jstar )=(T\jstar)(x),\qquad x\in X,$$
i.e., to find a fixed point of $T$ within ${\cal R}(X)$ (we can view $J^*=T\jstar$ as a generalized form of Bellman's equation). We also want to obtain a policy $\m^*\in{\cal M}$ such that $T_{\m^*}\jstar =T\jstar$. We assume that the control $u$ consists of the $m$ components $u_\ell$, $\ell=1,\ldots,m$, [cf.\ Eq.\ \controlstruct]. Note that since the state and control spaces are assumed finite, the control constraint set $U(x)$  and the set of policies ${\cal M}$ are also finite, so the minimum of various expressions over $U(x)$ or ${\cal M}$ is attained.

We will adopt throughout the following monotonicity  and contraction assumptions.

\xdef\assumptionmon{\assumptionn}\assumptionnum\show{myproposition}

\texshopbox{\assumption{\assumptionmon: (Monotonicity)} 
If $J,J'\in {\cal R}(X)$ and $J\le J'$, then
$$H(x,u,J)\le H(x,u,J'),\qquad \forall\ x\in X,\ u\in U(x).$$
}

\old{
Note that the monotonicity assumption implies the following properties, for all $J,J'\in {\cal R}(X)$ and $k=0,1,\ldots$:
$$J\le J'\qquad\implies\qquad T^kJ\le T^kJ',\qquad T_\m^kJ\le T_\m^kJ',\quad \forall\ \m\in{\cal M},$$
$$J\le TJ\qquad\implies\qquad T^kJ\le T^{k+1}J,\qquad T_\m^kJ\le T_\m^{k+1}J,\quad \forall\ \m\in{\cal M}.$$
Here $T^k$ and $T_\m^k$ denotes the $k$-fold composition of $T$ and $T_\m$, respectively.
}

For the contraction assumption, we introduce a function $v:X\mapsto\re$ with
$$v(x)>0,\qquad \forall\ x\in X.$$
We consider the weighted sup-norm
$$\|J\|=\max_{x\in X}{\big|J(x)\big|\over v(x)}$$
on  ${\cal R}(X)$, the space of real-valued functions $J$ on $X$.

\xdef\assumptionhc{\assumptionn}\assumptionnum\show{myproposition}

\texshopbox{\assumption{\assumptionhc: (Contraction)} 
For some $\a\in(0,1)$, we have
$$\|T_\m J-T_\m J'\|\le \a \|J-J'\|,\qquad \forall\ J,J'\in {\cal R}(X),\ \m\in{\cal M}.$$
}

The monotonicity and contraction  assumptions are satisfied in the $\a$-discounted finite-state MDP case \mdph, as well as other finite-state DP problems, such as stochastic shortest path problems in the special case where all policies are proper; see the books [BeT96], [Ber12], [Ber18] for an extensive discussion. In particular, for the $\a$-discounted MDP, $T_\m$ is a contraction with respect to the unweighted sup-norm with contraction modulus $\a$, whereas in the stochastic shortest path case, $T_\m$ is a contraction with respect to a weighted sup-norm with weights and contraction modulus that depend on the maximum expected time to reach the destination using proper policies (see [BeT96], Prop.\ 2.2). 

General abstract DP models under Assumptions \assumptionmon\ and \assumptionhc\ have been investigated in detail in the author's monograph [Ber18] [without assuming finiteness of $X$ and $U$, but with ${\cal R}(X)$ replaced by the set ${\cal B}(X)$ of all uniformly bounded functions over $X$, equipped with a weighted sup-norm]. The main results are that $T$ is a contraction mapping and has as unique fixed point the optimal cost function $\jstar$ (the equation $\jstar=T\jstar$ is Bellman's equation). Also $J_\m$ is the unique fixed point of $T_\m$. Moreover $\m$ is optimal if and only if $T_\m \jstar=T \jstar$ (or equivalently $T_\m J_\m=T J_\m$). Algorithmic results include the convergence of VI [i.e., $T^k J\to \jstar$ for all $J\in {\cal R}(X)$], and also convergence results for the PI algorithm \standardpi\ and some of its variations. We will be using these results in what follows in this paper, with the monograph [Ber18] as a general reference. For parts of our analysis, only the monotonicity and contraction Assumptions \assumptionmon\ and \assumptionhc\ are essential: the finiteness of the state and control spaces can be eliminated with minor mathematical proof modifications.

\section{Agent-by-Agent Value Iteration}

\pn The salient feature of the multiagent DP problem of this paper is that the control $u$ consists of $m$ components, 
$$u=(u_1,\ldots,u_m);$$
cf.\ Eq.\ \controlstruct. We will aim to develop a computationally efficient variant of the standard VI algorithm $J^{k+1}=TJ^k$, i.e.,
$$J^{k+1}(x)=\min_{(u_1,\ldots,u_m)\in U(x)}H(x,u_1,\ldots,u_m,J^k),\qquad x\in X.$$
Rather than simultaneous minimization over all the components $u_1,\ldots,u_m$, {\it our multiagent VI  algorithm involves sequential minimization of $H(x,u_1,\ldots,u_m,J^k)$ over a single component $u_\ell$, with the remaining components $u_{\ell'}$, $\ell'\ne \ell$, fixed at the values obtained through the preceding minimizations\/}. We maintain these control component values in a policy that is continually updated to incorporate the results of new minimizations.

Let $\m$ be a given policy that applies at $x$ the control
$$\m(x)=\big(\m_1(x),\ldots,\m_{m}(x)\big).$$
We define a constraint set for the $\ell$th control component $u_\ell$ that is  given by
$$U_{\ell,\m}(x)=\Big\{u_\ell\mid \big(\m_1(x),\ldots,\m_{\ell-1}(x),u_\ell,\m_{\ell+1}(x),\ldots,\m_m(x)\big)\in U(x)\Big\},\qquad \ell=1,\ldots,m.\xdef\compconstraint{\lab}\eqnum\show{oneo}$$
Note that since a policy $\m$ by definition satisfies the feasibility constraint
$$\big(\m_1(x),\ldots,\m_m(x)\big)\in U(x),\qquad x\in X,$$
the set $U_{\ell,\m}(x)$ contains $\m_\ell(x)$, so it is nonempty.
Note also that when $U(x)$ has the Cartesian product form $U_1(x)\times\cdots\times U_m(x)$, the set $U_{\ell,\m}(x)$ is simply equal to $U_\ell(x)$ for all $\m$.

Our algorithm generates a double sequence $\{J^k,\m^k\}$, starting from some pair $(J^0,\m^0)$: at the $k$th iteration, given $(J^k,\m^k)$, the algorithm obtains $(J^{k+1},\m^{k+1})$  after $m$ successive minimizations, one for each of the components $u_\ell$, $\ell=1,\ldots,m$. In particular, given the typical pair $(J,\m)$, our algorithm generates the next pair $(\tl J,\tl \m)$ as the last of a sequence of cost function-policy component pairs 
$$(\hat J_1,\hat \m_1),(\hat J_2,\hat \m_2),\ldots,(\hat J_m,\hat \m_m),\xdef\interpairs{\lab}\eqnum\show{oneo}$$
to be defined shortly, i.e., it sets
$$\tl J(x)=\hat J_m(x),\qquad \tl \m(x)=\big(\hat \m_{1}(x),\ldots,\hat \m_{m}(x)\big),\qquad x\in X.\xdef\finalresult{\lab}\eqnum\show{oneo}$$
The cost function-policy pairs \interpairs\ are obtained as follows:

For every $\ell=1,\ldots,m$, given $(\hat J_{\ell-1},\hat \m_{1},\ldots,\hat \m_{\ell-1})$, the algorithm generates $(\hat J_{\ell},\hat \m_{\ell})$ according to
$$\hat J_\ell(x)=\min_{u_\ell\in U_{\ell,(\hat \m_1,\ldots,\hat\m_{\ell-1},\m_\ell,\ldots,\m_{m})}(x)}H\big(x,\hat  \m_1(x),\ldots,\hat  \m_{\ell-1}(x),u_\ell,\m_{\ell+1}(x),\ldots,\m_m(x),\hat J_{\ell-1}\big),\qquad x\in X,\xdef\valiteragentbyagent{\lab}\eqnum\show{oneo}$$
$$\hat \m_\ell(x)\in\argmin_{u_\ell\in U_{\ell,(\hat \m_1,\ldots,\hat\m_{\ell-1},\m_\ell,\ldots,\m_{m})}(x)}H\big(x,\hat  \m_1(x),\ldots,\hat  \m_{\ell-1}(x),u_\ell,\m_{\ell+1}(x),\ldots,\m_m(x),\hat J_{\ell-1}\big),\qquad x\in X,\xdef\politeragentbyagenta{\lab}\eqnum\show{oneo}$$
where the constraint set in the two preceding minimizations,
$$U_{\ell,(\hat \m_1,\ldots,\hat\m_{\ell-1},\m_\ell,\ldots,\m_{m})}(x),$$
 is defined by Eq.\ \compconstraint; it is the set of $u_\ell$, which  are consistent (in terms of feasibilility) with the previously chosen components $\hat \m_1(x),\ldots,\hat\m_{\ell-1}(x)$ and the component choices  $\m_{\ell+1}(x),\ldots,\m_m(x)$ specified by the  policy $\m$. To start this process, only the initial function $\hat J_0$ is needed (in addition to $\m$), and it is given by
$$\hat J_0(x)=J(x),\qquad x\in X.\eqnum\show{oneo}$$

Note that each of the  minimizations \valiteragentbyagent\ is performed for every state $x\in X$, and that there may be multiple possible policies $\tl \m$ that can be generated by this process [cf.\ Eq.\ \finalresult], since the minimum in Eq.\ \politeragentbyagenta\ may not be uniquely attained. This set of policies is denoted by $\widetilde{\cal M}(J,\m)$.  Similarly, there may be multiple possible functions $\tl J$ that can be generated by this process [since the minimization \valiteragentbyagent\ is affected by the multiplicity of possible policies $\hat \m_1,\ldots,\hat\m_{\ell-1}$], and this set of functions is denoted by $\widetilde{\cal J}(J,\m)$.
In summary, our multiagent VI algorithm, starting from $(J^k,\m^k)$, generates $(J^{k+1}, \m^{k+1})$ according to
$$J^{k+1}\in\widetilde{\cal J}(J^k,\m^k),\qquad 
\m^{k+1}\in\widetilde{\cal M}(J^k,\m^k).\xdef\kthiteration{\lab}\eqnum\show{oneo}$$

\subsubsection{Optimistic and Asynchronous PI Algorithms}

\pn In the preceding algorithm \kthiteration, each iteration involves a policy improvement operation, i.e., an $m$-step minimization that cycles through all control components one-by-one. In Section 3, we will also consider an optimistic PI variant  where the $m$-step minimization is performed for only an infinite subset ${\cal K}\subset\{0,1,\ldots\}$ of the iterations, while for the complementary subset of iterations, $k\notin {\cal K}$, we use the (less expensive) standard policy evaluation update $J^{k+1}=T_{\m^k}J^k,$ and no policy update: 
$$J^{k+1}\in\widetilde{\cal J}(J^k,\m^k),\qquad 
\m^{k+1}\in\widetilde{\cal M}(J^k,\m^k),\qquad \forall\ k\in {\cal K},\xdef\optJiteration{\lab}\eqnum\show{oneo}$$
$$J^{k+1}=T_{\m^k}J^k,\qquad \m^{k+1}=\m^k,\qquad \forall\ k\notin {\cal K}.\xdef\optmuiteration{\lab}\eqnum\show{oneo}$$
We call the algorithm \optJiteration-\optmuiteration\  {\it multiagent optimistic PI\/}. It is a natural multiagent extension of the standard (single agent) optimistic PI algorithm \optpi, which is described in many sources for MDP and other problems, e.g., [Ber12], [Ber18], [Ber19]. Note that when ${\cal K}=\{0,1,\ldots\}$, the optimistic PI algorithm is the same as the multiagent VI algorithm \kthiteration. In cases where ${\cal K}$ is a ``small" subset of $\{0,1,\ldots\}$, the multiagent optimistic PI algorithm involves nearly exact policy evaluations and ``approaches" the multiagent PI algorithm proposed in our  earlier paper [Ber20].

In the preceding multiagent algorithms \kthiteration\ and \optJiteration-\optmuiteration, the iterations are performed simultaneously for all states $x\in X$. In Section 4, we will also consider an asynchronous distributed version of the multiagent optimistic PI algorithm 
 \optJiteration-\optmuiteration, whereby  iteration $k$ is performed for only a subset $X_k$ of the states. There is a requirement here is that each state $x$ belongs infinitely often to some subset $X_k$, so that there are infinitely many policy improvements at every state. This algorithm is well suited for distributed asynchronous computation, involving a partition of the state space into subsets, and with a processor assigned to each set of the partition.

\old{
$$J_1^k(x)=\min_{u_1\in U_1(x)}H\big(x,u_1,\m^k_2(x),\ldots,\m^k_m(x),J^k\big),\qquad x\in X,$$
$$\m^{k+1}_1(x)\in\argmin_{u_1\in U_1(x)}H\big(x,u_1,\m^k_2(x),\ldots,\m^k_m(x),J^k\big),\qquad x\in X,$$
$$\ldots\ \ \ \ldots$$
$$J_\ell^k(x)=\min_{u_\ell\in U_\ell(x)}H\big(x,\m^{k+1}_1(x),\ldots,\m^{k+1}_{\ell-1}(x),u_\ell,\m^k_{\ell+1}(x),\ldots,\m^k_m(x),J^k_{\ell-1}\big),\qquad x\in X,\xdef\valiteragentbyagent{\lab}\eqnum\show{oneo}$$
$$\m^{k+1}_\ell(x)\in\argmin_{u_\ell\in U_\ell(x)}H\big(x,\m^{k+1}_1(x),\ldots,\m^{k+1}_{\ell-1}(x),u_\ell,\m^k_{\ell+1}(x),\ldots,\m^k_m(x),J^k_{\ell-1}\big),\qquad x\in X,\eqnum\show{oneo}$$
$$\ldots\ \ \ \ldots$$
$$J_m^k(x)=\min_{u_m\in U_m(x)}H\big(x,\m^{k+1}_1(x),\ldots,\m^{k+1}_{m-1}(x),u_m,J_{m-1}^k\big),\qquad x\in X,$$
$$\m^{k+1}_m(x)\in\argmin_{u_m\in U_m(x)}H\big(x,\m^{k+1}_1(x),\ldots,\m^{k+1}_{m-1}(x),u_m,J_{m-1}^k\big),\qquad x\in X,$$
and finally,
$$J^{k+1}(x)=J_m^k(x),\qquad x\in X,\eqnum\show{oneo}$$
$$\m^{k+1}(x)=\big(\m_1^{k+1}(x),\ldots,\m_m^{k+1}(x)\big),\qquad x\in X.\eqnum\show{oneo}$$
}

\vskip-1pc

\section{Convergence to an Agent-by-Agent Optimal Policy}

\pn
We will prove that the  mutiagent VI algorithm \kthiteration\ converges to an agent-by-agent optimal policy, which we define as follows.

\xdef\definitionabgopt{\defn}\defnum\show{myproposition}

\texshopbox{\definition{\definitionabgopt: (Agent-by-Agent Optimality)}We say that a policy $\m=\{\m_1,\ldots,\m_m\}$ is {\it agent-by-agent optimal} if  for all $x\in X$ and $\ell=1,\ldots,m$, we have
$$H\big(x,\m_1(x),\ldots,\m_m(x),J_{\m}\big)=\min_{u_\ell\in U_{\ell,\m}(x)}H\big(x,\m_1(x),\ldots,\m_{\ell-1}(x),u_\ell,\m_{\ell+1}(x),\ldots,\m_m(x),J_{\m}\big),\xdef\agentbyagent{\lab}\eqnum\show{oneo}$$
where the constraint set $U_{\ell,\m}(x)$ is defined by Eq.\ \compconstraint.
} 

To interpret this definition, let a policy $\m=\{\m_1,\ldots,\m_m\}$ be given, and consider for every  $\ell=1,\ldots,m$ the single agent DP problem where for all $\ell'\ne \ell$ the $\ell'$th policy  components  are fixed at $\m_{\ell'}$, while the $\ell$th policy component is subject to optimization. The  definition \agentbyagent\ is the optimality condition for all the single agent problems; see [Ber18], Chapter 2 [Eq.\ \agentbyagent\ can be written as $T_{\m,\ell} J_\m=T_\ell J_\m$, where $T_\ell$ and $T_{\m,\ell}$ are the Bellman operators \tmap\ and \tmumap\ that correspond to the single agent problem involving agent $\ell$].  We can then conclude that $\m=\{\m_1,\ldots,\m_m\}$ is agent-by-agent optimal if each component $\m_\ell$ is optimal for the $\ell$th single agent problem, where it is assumed that the remaining policy components remain fixed; in other words by using $\m_\ell$, each agent $\ell$ acts optimally, assuming all other agents $\ell'\ne \ell$ continue to use the corresponding policy components $\m_{\ell'}$. 

Our definition of an agent-by-agent optimal policy is related to the notion of ``person-by-person" optimality from team theory, which has been studied primarily in the context of multiagent decision problems with nonclassical information patterns, whereby the agents may not share the information on which they base their decision. Thus team problems do not assume the shared information pattern that is characteristic of DP problems. For the origins of team theory and control with a nonclassical information pattern, we refer to Marschak [Mar55], Radner [Rad62], and Witsenhausen [Wit71a], [Wit71b], [Wit88]. For a sampling of subsequent works, we refer to the survey by Ho [Ho80], and the papers by Krainak, Speyer, and Marcus [KLM82a], [KLM82b], de Waal and van Schuppen [WaS00]. For more  recent works, see Nayyar, Mahajan, and Teneketzis [NMT13], Nayyar and Teneketzis [NaT19],  Li et al.\ [LTZ19], Gupta [Gu20], the book by Zoppoli, Parisini, Baglietto, and Sanguineti [ZPB19], and the references quoted there. 

Note that an (overall) optimal policy is agent-by-agent optimal, but the reverse may not be true. This is similar to properties of person-by-person optimal solutions in team theory. It is also similar to what may happen in coordinate descent methods for multivariable optimization, where it is possible (in the absence of favorable assumptions) to stop at a nonoptimal point where no progress can be made along any one coordinate; some  examples involving a Cartesian product constraint set of the form \cartprod\ are given in the paper [Ber20]. 

While an agent-by-agent optimal policy may be either optimal or adequate for practical purposes, it may offer no guarantees of quality. For a simple example, let $U(x)$ be the intersection of a Cartesian product of finite subsets $U_\ell(x)$ of the real line and the unit simplex:
$$U(x)=\big\{U_1(x)\times\cdots\times U_m(x)\big\}\cap \{u\mid u_1+\cdots+u_m=1\}.$$
Then it can be seen that  the constraint set $U_{\ell,\m}(x)$ consists of just the single point $\m_\ell(x)$, so that {\it all} feasible policies are agent-by-agent optimal. This is due to the extreme coupling of the control components through the simplex constraint. It would not happen if the constraint set was just a Cartesian product $U_1(x)\times\cdots\times U_m(x)$, in which case $U_{\ell,\m}(x)=U_\ell(x)$ for all $\ell$. Nonetheless, one should be aware that the method of partitioning of the control into components may seriously impact the effectiveness of  our multiagent VI algorithm through the creation of spurious agent-by-agent optimal policies.

\xdef\assumptionunique{\assumptionn}\assumptionnum\show{myproposition}

We will now prove our main convergence result, under the following assumption, which is reminiscent of strict convexity assumptions in the analysis  of coordinate descent methods (see e.g., [Ber16], Section 3.7). While we do not have a concrete counterexample, we speculate based on experience with coordinate descent methods, that the assumption cannot be easily dispensed with.

\texshopbox{\assumption{\assumptionunique: (Uniqueness Property)}The cost functions of district policies are distinct, i.e., for any two policies $\m$ and $\m'$
$$\m\ne\m'\qquad\implies\qquad J_\m\ne J_{\m'}.$$
}

Our convergence result also assumes that the initial condition $(J^0,\m^0)$ satisfies
$$T_{\m^0}J^0\le J^0.\xdef\initcond{\lab}\eqnum\show{oneo}$$
This assumption is unnecessary for the $\a$-discounted MDP where $T$ and $T_\m$ are given by Eqs.\ \tmap\ and \tmumap. The reason is that if we replace $J^0$ by a function $\ol J^0$ obtained by shifting $J^0$ by a constant $c$ [i.e., replace $J^0(x)$ by $J^0(x)+c$ for all $x$], we will have 
$$(T_{\m^0}\ol J^0)(x)=(T_{\m^0}J^0)(x)+\a c\le J^0(x)+c=\ol J^0(x),$$
provided $c$ is large enough, thereby satisfying the assumption \initcond. At the same time, it can be seen that by replacing  $J^0$ with $\ol J^0$ the generated policies will not be affected, while the generated iterates $J^k$ will just be shifted by an appropriate constant. Thus for discounted MDP the assumption \initcond\ is unnecessary for the following convergence result, since the same sequence of policies will be obtained whether we use $J^0$ or $\ol J^0$. 

For other types of problems the assumption \initcond\ is needed. However, thanks to the contraction property of Assumption \assumptionhc, it can be typically satisfied by adding to $J^0(x)$ a sufficiently large constant $c$ for all $x$. In particular, any function $\ol J$ that satisfies
$$\a\|\ol J-J_\m\| \le {\ol J(x)\over v(x)}-{J_\m(x)\over v(x)},\qquad \forall\ x\in X,\ \m\in {\cal M},\xdef\oljcond{\lab}\eqnum\show{oneo}$$ 
(for example a sufficiently large constant function) also satisfies the condition \initcond. To see this, note that for $J\le \ol J$ and $x\in X$,
$${(T_\m J)(x)\over v(x)}\le {(T_\m\ol J)(x)\over v(x)}\le {(T_\m J_\m)(x)\over v(x)}+\a\|\ol J-J_\m\|= {J_\m(x)\over v(x)}+\a\|\ol J-J_\m\|\le {\ol J(x)\over v(x)},$$
where the first inequality follows from the monotonicity of $H$, the second inequality follows by applying the contraction property with $J=\ol J$, $J'=J_\m$, and the third inequality is Eq.\ \oljcond. Thus, for $\ol J$ satisfying Eq.\ \oljcond, we have $T_\m\ol J\le \ol J$ for all $\m\in {\cal M}$.

\xdef\propmultiagentviter{\propn}\propnum\show{myproposition}

\texshopbox{\proposition{\propmultiagentviter: (VI Convergence to an Agent-by-Agent Optimal Policy)}Let Assumptions \assumptionmon, \assumptionhc, and \assumptionunique\ hold, and assume further that the state and control spaces $X$ and $U$ are finite, and that the initial pair $(J^0,\m^0)$ satisfies Eq.\ \initcond. Let $\{J^k,\m^k\}$ be a sequence generated by the agent-by-agent VI algorithm \kthiteration. Then there is an agent-by-agent optimal policy $\bar \m$  and an index $\ol k$ such that for all $k\ge \ol k$, we have $\m^k=\bar \m$, and
$$\|J^{k+1}-J_{\bar \m}\|\le \a\|J^{k}-J_{\bar \m}\|,\xdef\convrate{\lab}\eqnum\show{oneo}$$
while the sequence $\{J^k\}$ converges to $J_{\bar \m}$.}

\proof The critical step of the proof is to show that for all $(J,\m)$ with
$$T_\m J\le J,$$
and all $(\tl J,\tl\m)$ with $\tl J\in \widetilde {\cal J}({J,\m})$ and $\tl\m\in \widetilde {\cal M}({J,\m})$ [cf.\ Eq.\ \kthiteration],
the following monotone decrease inequality holds 
$$T_{\tl \m}\tl J\le \tl J=\hat J_{m}\le\hat J_{m-1}\le \hat J_{m-2}\le \cdots\le \hat J_1\le T_{\m}J\le J,\xdef\monotonedecrease{\lab}\eqnum\show{oneo}$$
where for all $\ell=1,\ldots,m,$ and  $x\in X$,
$$\eqalign{\hat J_\ell(x)&=\big(T_{(\hat \m_1,\ldots,\hat \m_\ell,\m_{\ell+1},\ldots,\m_m)}\hat J_{\ell-1}\big)(x)\cr
&=\min_{u_\ell\in U_{\ell,(\hat \m_1,\ldots,\hat\m_{\ell-1},\m_\ell,\ldots,\m_{m})}(x)}H\big(x,\hat \m_1(x),\ldots,\hat \m_{\ell-1}(x),u_\ell,\m_{\ell+1}(x),\ldots,\m_m(x),\hat J_{\ell-1}\big),\cr}\eqnum\show{oneo}$$
with $\hat J_0=J$  [cf.\ Eq.\ \valiteragentbyagent], and 
$$\hat \m_\ell(x)\in\argmin_{u_\ell\in U_{\ell,(\hat \m_1,\ldots,\hat\m_{\ell-1},\m_\ell,\ldots,\m_{m})}(x)}H\big(x,\hat  \m_1(x),\ldots,\hat  \m_{\ell-1}(x),u_\ell,\m_{\ell+1}(x),\ldots,\m_m(x),\hat J_{\ell-1}\big),$$
[cf.\ Eq.\ \politeragentbyagenta]. Indeed the relation \monotonedecrease\ is proved starting from the right side, which is the assumption \initcond, and by using the definition of the algorithm, and the monotonicity Assumption \assumptionmon\ to prove first that $\hat J_1\le T_{\m}J\le J$, and then by proceeding sequentially to the inequality $\hat J_{m}\le\hat J_{m-1}$.  
 In particular, at the typical step, assuming that $\hat J_{\ell-1}\le \hat J_{\ell-2}$, we show that $\hat J_{\ell}\le \hat J_{\ell-1}$ by writing
$$\eqalign{\hat J_{\ell-1}(x)&=H\big(x,\hat \m_1(x),\ldots,\hat \m_{\ell-1}(x),\m_\ell(x),\m_{\ell+1}(x),\ldots,\m_m(x),\hat J_{\ell-2}\big),\cr
&\ge H\big(x,\hat \m_1(x),\ldots,\hat \m_{\ell-1}(x),\m_\ell(x),\m_{\ell+1}(x),\ldots,\m_m(x),\hat J_{\ell-1}\big),\cr
&\ge \min_{u_\ell\in U_{\ell,(\hat \m_1,\ldots,\hat\m_{\ell-1},\m_\ell,\ldots,\m_{m})}(x)}H\big(x,\hat \m_1(x),\ldots,\hat \m_{\ell-1}(x),u_\ell,\m_{\ell+1}(x),\ldots,\m_m(x),\hat J_{\ell-1}\big),\cr
&=\hat J_{\ell}(x),\cr}$$
where the first inequality follows by using the monotonicity Assumption \assumptionmon\ and the hypothesis $\hat J_{\ell-1}\le \hat J_{\ell-2}$.
Finally, by applying $T_{\tl \m}$ to the relation $\hat J_{m}\le\hat J_{m-1}$ to obtain $T_{\tl \m}\hat J_{m}\le T_{\tl \m}\hat J_{m-1}$, and by using the facts
$\tl J=\hat J_{m}=T_{\tl \m}\hat J_{m-1}$, we obtain the leftmost relation $T_{\tl \m}\tl J\le \tl J=\hat J_{m}$ in Eq.\monotonedecrease. 
(Note that the contraction assumption is not needed for the preceding argument, and this is useful for applying this line of proof in other DP problem contexts.)

From Eq.\ \monotonedecrease, we see that the sequence of functions $J^k$ converges monotonically to some function $\ol J$, and the same is true for all the sequences of intermediate functions $J^k_1,\ldots,J^k_{m-1}$. 
For each $\ell$, let the policies 
$$(\m_1^{k+1},\ldots,\m_\ell^{k+1},\m_{\ell+1}^k,\ldots,\m_m^k)$$
 be equal to some policy $\bar\m[\ell]=\big(\bar\m_1[\ell],\ldots,\bar\m_m[\ell]\big)$ infinitely often  (such a policy exists since the set of policies is finite). Then we will have for all $x\in X$ and $\ell=1,\ldots,m$,
$$J^k_\ell(x)=H(x,\bar\m[\ell](x),J^k_{\ell-1})=\min_{u_\ell\in U_{\ell,\bar\m[\ell]}(x)}H\big(x,\bar\m_1[\ell](x),\ldots,\bar\m_{\ell-1}[\ell](x),u_\ell,\bar\m_{\ell+1}[\ell](x),\ldots,\bar\m_{m}[\ell](x),J^k_{\ell-1}\big),$$
infinitely often.
By taking limit as $k\to\infty$ and using the continuity of $H(x,u,\cdot)$ (which is implied by the contraction property of $T_{\bar\m[\ell]}$), we have
$$\ol J(x)=H\big(x,\bar\m[\ell](x),\ol J\big)=(T_{\bar\m[\ell]}\ol J)(x),\qquad \ell=1,\ldots,m,\ x\in X,\xdef\firstcond{\lab}\eqnum\show{oneo}$$
as well as
$$\ol J(x)=\min_{u_\ell\in U_{\ell,\bar\m[\ell]}(x)}H\big(x,\bar\m_1[\ell](x),\ldots,\bar\m_{\ell-1}[\ell](x),u_\ell,\bar\m_{\ell+1}[\ell](x),\ldots,\bar\m_{m}[\ell](x),\ol J\big),\qquad \ell=1,\ldots,m,\ x\in X.\xdef\secondcond{\lab}\eqnum\show{oneo}$$
Equation \firstcond\ and the contraction property of $T_{\bar\m[\ell]}$ imply that $\ol J$ is equal to the cost functions $J_{\bar\m[\ell]}$ of all of the $m$ policies $\bar\m[\ell]$, $\ell=1,\ldots,m$. In view of the uniqueness Assumption \assumptionunique, this implies that all the policies $\bar\m[\ell]$, $\ell=1,\ldots,m$, are equal to some policy $\bar\m$, which has cost function $\ol J$, and in view of Eq.\ \secondcond, satisfies
$$\ol J(x)=\min_{u_\ell\in U_{\ell,\bar\m}(x)}H\big(x,\bar\m_1(x),\ldots,\bar\m_{\ell-1}(x),u_\ell,\bar\m_{\ell+1}(x),\ldots,\bar\m_{m}(x),\ol J\big),\qquad \ell=1,\ldots,m.\xdef\abamin{\lab}\eqnum\show{oneo}$$
It follows that $\bar\m$ is agent-by-agent optimal. 

Finally, the preceding argument shows that $\ol J$ is the cost function of every policy that is repeated infinitely often. Thus the uniqueness Assumption \assumptionunique\ implies that $\bar\m$ is the only policy that is repeated infinitely often. Since there are finitely many policies, it follows that $\m^k=\bar\m$ for all $k$ after some index. Hence from the definition of the algorithm, the sequence $\{J^k\}$ satisfies 
$J^{k+1}=T_{\bar\m}J^k$
for all $k$ after some index, which in view of the contraction Assumption \assumptionhc, implies Eq.\ \convrate.
\qed

Note that the preceding proposition does not guarantee convergence  to the optimal policy (which is unique by Assumption  \assumptionunique). In particular, if our algorithm is started at a pair $(J_{\bar\m},\bar\m)$, where $\bar\m$ is an agent-by-agent optimal policy, it will not move from $\bar\m$ [in fact it can be shown that this will happen even if the algorithm is started at a pair $(J^0,\bar\m)$, where $J^0$ is sufficiently close to $J_{\bar\m}$]. Thus every agent-by-agent optimal policy behaves like a ``local minimum," with its own ``region of attraction," and is a potential convergence limit of our algorithm. The limit will depend on the starting pair, as well as the order in which the agents select their components. The algorithm guarantees convergence  to the optimal policy only under additional assumptions that guarantee that there are no additional agent-by-agent optimal policies.  We postpone a discussion of this issue for Section 5.

\subsubsection{Ensuring Convergence to an Optimal Policy with Randomization Schemes}

\pn Another possibility to enhance the convergence properties of the algorithm, and ensure convergence to an optimal policy, is to enlarge the constraint sets
$$U_{\ell,(\hat \m_1,\ldots,\hat\m_{\ell-1},\m_\ell,\ldots,\m_{m})}(x)$$
in Eq.\ \valiteragentbyagent, to allow minimization over subsets of multiple control components. These subsets may be selected with some form of randomization: at some iterations minimize over a single control component as in iteration \valiteragentbyagent-\politeragentbyagenta, while at some other randomly chosen iterations minimize over multiple or even all control components. Schemes of this type  have been considered for the purpose of enhancing the convergence properties of asynchronous PI; see [Ber18], Section 2.5.3. Randomization over sets of multiple control components can also be used in the context of the optimistic agent-by-agent PI methods of the next section, and they can similarly enhance their convergence properties. 

We will not consider randomized control component selection schemes in this paper. Their analysis is  similar to the one of [Ber18], Section 2.5.3, their implementation is likely problem-dependent, and their practical performance is an interesting subject for further research.
 Their principal drawback is that simultaneous minimization over multiple control components can be very costly (depending on the number of components involved), even if it used in only a small proportion of the total number of iterations. 
\vskip-1pc

\section{Agent-by-Agent Optimistic Policy Iteration}

\pn Let us now consider an optimistic PI variant  where we introduce an infinite subset ${\cal K}\subset\{0,1,\ldots\}$ of the iterations, and the complementary subset of iterations $k\notin {\cal K}$. For the latter subset, we use the (less expensive) standard policy evaluation update $J^{k+1}=T_{\m^k}J^k,$ and no policy update: 
$$J^{k+1}\in\widetilde{\cal J}(J^k,\m^k),\qquad 
\m^{k+1}\in\widetilde{\cal M}(J^k,\m^k),\qquad \forall\ k\in {\cal K},\xdef\optJiteration{\lab}\eqnum\show{oneo}$$
$$J^{k+1}=T_{\m^k}J^k,\qquad \m^{k+1}=\m^k,\qquad \forall\ k\notin {\cal K}.\xdef\optmuiteration{\lab}\eqnum\show{oneo}$$
We have the following convergence result:

\xdef\propmultiagentoptpi{\propn}\propnum\show{myproposition}

\texshopbox{\proposition{\propmultiagentoptpi: (Optimistic PI Convergence to an Agent-by-Agent Optimal Policy)}Let the assumptions of  Prop.\ \propmultiagentviter\ hold, and let $\{J^k,\m^k\}$ be a sequence generated by the optimistic agent-by-agent PI algorithm \optJiteration-\optmuiteration. Then there is an index $\ol k$ such that for all $k\ge \ol k$, we will have $\m^k=\bar \m$, where  $\bar \m$  is an agent-by-agent optimal policy, while the sequence $\{J^k\}$ will converge to $J_{\bar \m}$.}

\proof The proof is essentially identical to the one of Prop.\ \propmultiagentviter. \qed

The algorithm admits also a distributed implementation, whereby the iteration \optJiteration-\optmuiteration\ is executed at the subset of times $k\in {\cal K}$ only for a subset $X_k$ of the states, while for the remaining states $x\notin{\cal K}$ the values of $J^{k+1}(x)$ and $\m^{k+1}(x)$ remain unchanged:
$$J^{k+1}(x)=J^{k}(x),\qquad \m^{k+1}(x)=\m^k(x),\qquad \forall\ x\notin X_k,\ k\in {\cal K}.\eqnum\show{oneo}$$
In addition to the set ${\cal K}$ being infinite, there is a requirement here is that each state $x$ belongs infinitely often to some subset $X_k$, so that there are infinitely many policy improvements at every state. Algorithms of this type have been proposed in the book [BeT96], Section 2.2.3, and in [Ber12].
The convergence proof of Prop.\ \propmultiagentviter\ still goes through; see also the proof of Prop.\ 2.5 of [BeT96]. Note, however, that for this type of algorithm to be provably convergent, $(J^0,\m^0)$ must satisfy the condition  $T_{\m^0}J^0\le J^0$ [cf.\ Eq.\ \initcond] even for discounted MDP, as demonstrated with counterexamples  by Williams and Baird [WiB93] (see also [Ber10]).

In a more complex version of the algorithm, the information on the cost function iterates at each iteration is allowed to be out-of-date, while modifications are introduced to eliminate the need for the initial condition assumption of Eq.\ \initcond. Distributed asynchronous PI algorithms of this type have been proposed and analyzed in the paper by Bertsekas and Yu [BeY10]  [see also [BeY12], [YuB13], and the books [Ber12] (Section 2.6), and [Ber18] (Section 2.6)]. See also the randomized optimistic PI algorithms of [Ber18] (Section 2.5.3). Multiagent versions of such algorithms are a subject for further research.

\vskip-1pc

\section{Conditions for Obtaining an Optimal Policy}
\vskip-0.5pc

\pn We proved earlier that our multiagent VI algorithm will find an agent-by-agent optimal policy under our assumptions of Prop.\ \propmultiagentviter, but this policy need not be optimal. We will now discuss approaches that can be used to show that the policy obtained is optimal, under the same or alternative assumptions. One possibility is to impose conditions under which every agent-by-agent optimal policy is optimal. To this end we introduce the following definition.

\xdef\definitioncompbycomp{\defn}\defnum\show{myproposition}

\texshopbox{\definition{\definitioncompbycomp: (Component-by-Component Minimum)} For a state $x$ and a function $J\in {\cal R}(X)$ we say that a control $\ol u=(\ol u_1,\ldots,\ol u_m)\in U(x)$ is a {\it component-by-component minimum of $H$ at $(x,J)$} if
$$\ol u_\ell\in\arg\min_{u_\ell\in \ol U_{\ell,\bar u}(x)} H(x,\ol u_1,\ldots,\ol u_{\ell-1},u_\ell,\ol u_{\ell+1},\ldots,\ol u_m,J),\qquad \ell=1,\ldots,m.\eqnum\show{oneo}$$
where the sets $\ol U_{\ell,\bar u}(x)$ are defined by
$$\ol U_{\ell,\bar u}(x)=\big\{u_\ell\mid (\ol u_1,\ldots,\ol u_{\ell-1},u_\ell,\ol u_{\ell+1},\ldots,\ol u_m)\in U(x)\big\},\qquad \ell=1,\ldots,m.$$
} 

Note that from the definition of agent-by-agent optimality, we have that $\bar \m$ is agent-by-agent optimal if for every $x\in X$, the control $\bar\m(x)$ is a component-by-component minimum of $H$ at $(x,J_{\bar\m})$.
We have the following proposition.

\xdef\propagentbyagen{\propn}\propnum\show{myproposition}

\texshopbox{\proposition{\propagentbyagen: (Agent-by-Agent Optimality Criterion)}
\nitem{(a)} Assume that for every state $x\in X$ and policy $\bar \m$ such that $\bar \m(x)$ is a component-by-component minimum of $H$ at $(x,J_{\bar\m})$, the control $\bar\m(x)$ minimizes $H(x,u,J_{\bar\m})$ over $u\in U(x)$. Then every agent-by-agent optimal policy is optimal.
\nitem{(b)}  If an agent-by-agent optimal policy $\bar \m$ is optimal, then for all $x\in X$ the control $\bar\m(x)$ is a component-by-component minimum of $H$ at all $(x,J_{\bar\m})$.
}

\proof (a) Let $\bar\m$ be agent-by-agent optimal. Then from the definition of agent-by-agent optimality, we have that for all $x\in X$, $\bar\m(x)$ is a component-by-component minimum of $H$ at $(x,J_{\bar\m})$. By our assumption, this implies that for all $x\in X$, $\bar\m(x)$ minimizes $H(x,u,J_{\bar\m})$ over $u\in U(x)$, or  $T_{\bar\m}J_{\bar\m}=TJ_{\bar\m}$. From general properties of contractive abstract DP models (cf.\ [Ber18], Chapter 2), we also have $T_{\bar\m}J_{\bar\m}=J_{\bar\m}$. Hence $T_{\bar\m}J_{\bar\m}=TJ_{\bar\m}$, which implies that $J_{\bar\m}=\jstar$ (cf.\ [Ber18], Chapter 2), so $\bar\m$ is optimal.
\smskip
\pn (b) Optimality of $\bar\m$ implies that $T_{\bar\m}J_{\bar\m}=TJ_{\bar\m}$  (cf.\ [Ber18], Chapter 2). This implies that for all $x\in X$ the control $\bar\m(x)$ is a component-by-component minimum of $H$ at all $(x,J_{\bar\m})$. \qed
\smskip

In view of Prop.\ \propagentbyagen(a), an important issue is to delineate  sufficient conditions that guarantee that component-by-component minima of $H$ at $(x,J_\m)$ minimize $H(x,u,J_\m)$ over $u\in U(x)$. Somewhat similar questions have been addressed in two related contexts:
\nitem{(a)} Team theory in connection with the notion of person-by-person optimality mentioned earlier.
\nitem{(b)} The theory of convergence of coordinate descent methods in nonlinear optimization.
\smskip

In the theory of teams and other related decentralized control problem formulations,  the most prominent analytical issues arise when the team members select control components based on different information. By contrast in our framework the agents choose actions based on shared information, namely the current state $x_k$ of the system. Because of this fundamental structural assumption, DP algorithms such as VI and PI apply to our problem, but do not apply to team problems with nonclassical information pattens, exemplified by the famous counterexample of Witsenhausen [Wit68]. 

In the theory of coordinate descent methods, the result most related to our context is that if a function $F(y_1,\ldots,y_m)$ of $m$ vectors  $y_1,\ldots,y_m$ is strictly convex and differentiable over the Cartesian product $Y_1\times\cdots\times Y_m$ of closed convex sets $Y_1,\ldots,Y_m$, then a vector $\bar y=(\bar y_1,\ldots,\bar y_m)$ is a global minimum of $F$ over $Y_1\times\cdots\times Y_m$ if and only if it has the component-by-component minimization property
$$\bar y_\ell\in\arg\min_{y_\ell\in Y_\ell}F(\bar y_1,\ldots,\bar y_{\ell-1},y_\ell,\bar y_{\ell+1},\ldots,\bar y_m),\qquad\hbox{for all }\ell=1,\ldots,m.\xdef\cornercond{\lab}\eqnum\show{oneo}$$
Thus when $F$ is strictly convex and differentiable, the block coordinate descent method cannot get trapped into a solution that is component-by-component minimum but is not a global minimum [this is not true, however, if $F$ is strictly convex but nondifferentiable, since the condition \cornercond\ may hold at vectors $\bar y$ that are not global minima, and at which $F$ is nondifferentiable]. Some related results are known for the case where the sets $Y_\ell$ are discrete, under assumptions that can be viewed as discrete space substitutes for strict convexity; see e.g., de Waal and van Schuppen [WaS00], and Bauso and Pesenti [BaP08], [BaP12]. 

While the coordinate descent and the team theory results provide some analytical guidance, they do not apply directly to the DP context of this paper. The reason is that the mapping $H$ involves the functions $J_\m$, whose properties have to be verified through analysis. 
We leave this line of investigation as a subject for further research, and we outline another analytical approach, which assumes continuous state and control spaces $X$ and $U$, and is based on strict convexity and differentiability assumptions. 

\subsubsection{Continuous Spaces, Strict Convexity, and Differentiability}

\pn Let us remove the assumption that the state and control spaces $X$ and $U$ are finite, while continuing to assume that the control has $m$ components, $u=(u_1,\ldots,u_m)$ that are constrained by $u\in U(x)$ for all $x\in X$. We continue to adopt the monotonicity and contraction Assumptions \assumptionmon\ and \assumptionhc, with the modification that ${\cal R}(X)$ is replaced by the space ${\cal B}(X)$ of bounded functions over $X$, with respect to a weighted sup-norm. Moreover, we assume that the various minima over control components in the definition of the algorithms are attained. Models of this type have been analyzed extensively in the monograph [Ber18] (Chapter 2), to which we refer for a detailed discussion.
The definitions of agent-by-agent optimality and component-by-component minimum  carry over without change to the continuous spaces setting, and so does the associated agent-by-agent optimality criterion [cf.\ Prop.\ \propagentbyagen(a)]. Furthermore, the key inequality \convrate\ for the proof of the convergence result of Prop.\ \propmultiagentviter\ goes through, under the condition $T_{\m^0}J^0\le J^0$ [cf.\ Eq.\ \initcond]. As a result, the proof of monotonic decrease of the sequence $\{J^k\}$ to some function $\ol J$ goes through as well.

In conclusion, without assuming finiteness of the state and control spaces $X$ and $U$, our algorithm, under the monotonicity  and contraction Assumptions \assumptionmon\ and \assumptionhc, and the condition $T_{\m^0}J^0\le J^0$ [cf.\ Eq.\ \initcond], converges monotonically to some $\ol J$, which can be seen to be pointwise bounded below by the optimal cost function $\jstar$, which belongs to ${\cal B}(X)$, so that $\ol J\in {\cal B}(X)$. Further conditions, involving  strict convexity and differentiability, need to be imposed to guarantee that $\ol J=\jstar$, that $\jstar$ is convex and differentiable, and that an optimal policy can be obtained. A stochastic optimal control model, involving a linear system, a convex cost per stage, and convex state and control constraints,  was formulated and analyzed in 1973 by the author [Ber73], and  is well suited for this purpose. We leave further analysis along this line as a subject for further research.

\old{
TO BE CONTINUED FROM HERE WITH A CONVEX STOCHASTIC OPTIMAL CONTROL PROBLEM
Let us now return to our  DP problem under the monotonicity, contraction, and uniqueness Assumptions  \assumptionmon, \assumptionhc, and \assumptionunique. Since under these assumptions there is a unique optimal policy, if we can show that there is a unique agent-by-agent optimal policy, this policy must be optimal. Alternatively, we can try to verify the condition of Prop.\ \propagentbyagen(a). If we could argue that for each $(x,J_\m)$, the function $H(x,u,J_\m)$, viewed as a function of $u$, is strictly convex and differentiable, there would be a unique agent-by-agent optimal policy. This may be useful in a more general formulation of our problem that allows continuous control spaces, but not under our assumption of finiteness of the control space. 
The approach of using strict convexity and differentiability can be refined by restricting attention to a subset of initial conditions for our algorithm. Assume that we have a subset of strictly convex and differentiable functions ${\cal F}\subset {\cal R}(X)$ with the property that
if the multiagent VI algorithm \kthiteration\ is started within ${\cal F}$, it remains within ${\cal F}$, i.e., if $\{J^k\}$ is a generated sequence by the algorithm, and $J^0\in {\cal F}$, then $\{J^k\}\subset {\cal F}$, and 
furthermore ${\cal F}$ is closed in the sense that it contains the pointwise limits of all convergent sequences that belong to ${\cal F}$
Then it follows that if $J^0\in {\cal F}$, the sequence $\{J^k\}$ generated by the algorithm is contained in ${\cal F}$, and so is the limit function $\ol J$. Moreover, an examination of the proof of Prop.\ \propmultiagentviter\ shows that the component-by-component minimization equation
$$\ol J(x)=\min_{u_\ell\in U_{\ell,\bar\m}(x)}H\big(x,\bar\m_1(x),\ldots,\bar\m_{\ell-1}(x),u_\ell,\bar\m_{\ell+1}(x),\ldots,\bar\m_{m}(x),\ol J\big),\qquad \ell=1,\ldots,m,$$
holds for the policy $\bar \m$ that is obtained after a finite number of iterations [cf.\ Eq.\ \abamin]. It then follows that $\bar \m$ is an optimal policy [cf.\ the proof of Prop.\ \propagentbyagen(a)].
}

\old{DOES NOT WORK
\subsubsection{An Alternative Approach for Showing that the Limit Policy is Optimal}
\pn Let us consider an alternative approach for showing that the policy $\bar \m$ obtained by our algorithm is optimal, thereby strengthening the convergence result of Prop.\ \propmultiagentviter. Assume that we have a set of functions ${\cal F}\subset {\cal R}(X)$ with the following properties:
\nitem{(a)} If the multiagent VI algorithm \kthiteration\ is started within ${\cal F}$, it remains within ${\cal F}$, i.e., if $\{J^k\}$ is a generated sequence by the algorithm, and $J^0\in {\cal F}$, then $\{J^k\}\subset {\cal F}$.
\nitem{(b)} ${\cal F}$ is closed in the sense that it contains the limits of all convergent sequences that belong to ${\cal F}$. 
\nitem{(b)} ${\cal F}$ has the property that for any $x\in X$ and $J\in{\cal F}$, every component-by-component minimum $\ol u=(\ol u_1,\ldots,\ol u_m)\in U(x)$ of $H$ at $(x,J)$ minimizes $H(x,u,J)$ over $u\in U(x)$. 
\smskip
From properties (a) and (b) above, it follows that if $J^0\in {\cal F}$, the sequence $\{J^k\}$ generated by the multiagent VI algorithm \kthiteration\ is contained in ${\cal F}$, and so is the limit function $\ol J$. Moreover, an examination of the proof of Prop.\ \propmultiagentviter\ shows that the component-by-component minimization equation
$$\ol J(x)=\min_{u_\ell\in U_{\ell,\bar\m}(x)}H\big(x,\bar\m_1(x),\ldots,\bar\m_{\ell-1}(x),u_\ell,\bar\m_{\ell+1}(x),\ldots,\bar\m_{m}(x),\ol J\big),\qquad \ell=1,\ldots,m,$$
holds for the policy $\bar \m$ that is obtained after a finite number of iterations [cf.\ Eq.\ \abamin]. From property (c) above, it then follows that $\bar \m$ is an optimal policy [cf.\ the proof of Prop.\ \propagentbyagen(a)].
In conclusion, if there is a suitable set ${\cal F}$ for which properties (a)-(c) can be verified, and the initial condition $J^0$ belongs to ${\cal F}$, the multiagent VI algorithm \kthiteration\ yields an optimal policy under the assumptions of Prop.\ \propmultiagentviter. This is true even if there exist some additional (suboptimal) agent-by-agent policies.
The preceding approach leaves the burden of defining the set ${\cal F}$ and verifying of properties (a)-(c) above, to a problem-dependent analysis. For some problems with favorable structure the approach may be fruitful, while for others it may not. We provide an example of a problem with favorable structure. 
\xdef\examplemdpmon{\exampl}\examplnum\show{myexample}
\beginexample{\examplemdpmon\ (Monotone $\a$-Discounted MDP with Cartesian Product Control Structure)}Consider an $\a$-discounted MDP where $H$ has the form
$$H(x,u,J)=\sum_{y=1}^np_{xy}(u)\big(g(x,u,y)+\a J(y)\big),\qquad x=1,\ldots,n;$$
cf.\ Eq.\ \mdph, and we assume the following:
\nitem{(1)} $U=U(x)=U_1\times\cdots\times U_m$.
\nitem{(2)} ${\cal F}$ consists of all $J\in {\cal R}(X)$ such that for any two controls $u,\bar u\in U$, 
$$u\le \bar u\qquad\implies\qquad  \sum_{y=1}^np_{xy}(u)J(y)\le \sum_{y=1}^np_{xy}(\bar u)J(y),\quad \hbox{for all }x\in X.$$
\nitem{(3)} For any two controls $u,\bar u\in U$,
$$u\le \bar u\qquad\implies\qquad \sum_{y=1}^np_{xy}(u)g(x,u,y)\le \sum_{y=1}^np_{xy}(\bar u)g(x,\bar u,y),\quad \hbox{for all }x\in X.$$
\nitem{(4)} For any two states $x,\bar x\in X$, 
$$x\le \bar x\qquad\implies\qquad \sum_{y=1}^np_{\bar xy}(u)g(x,u,y)\le \sum_{y=1}^n p_{\bar xy}(u)g(\bar x,u,y),\quad \hbox{for all }u\in U.$$
\smskip
Then it can be shown that properties (a) and (b) above are satisfied and the policy $\bar \m$ obtained by our multiagent VI algorithm \kthiteration\ is optimal. Indeed property (a) clearly holds, so we focus on proving property (b), namely that
for any $x\in X$ and $J\in{\cal F}$, every component-by-component minimum $\ol u=(\ol u_1,\ldots,\ol u_m)\in U(x)$ of $H$ at $(x,J)$ minimizes $H(x,u,J)$ over $u\in U(x)$.
\endexample
}

\old{DOES NOT WORK
In what follows, we discuss a special case where the component-by-component minimum property can be verified, so that Prop.\ \propagentbyagen(a) applies and shows that there is a unique agent-by-agent optimal policy, which must be optimal. 
\subsubsection{Monotone Discrete DP Problems}
\pn We will introduce a dynamic version of a static multivariable discrete optimization problem that was formulated and analyzed in the paper by Wu and Bertsekas [WuB01] (see also [Ber16], Exercise 6.5.4).
Consider an $\a$-discounted  deterministic optimal control problem [$\a\in(0,1)$], where $H$ is given by
$$H(x,u,J)=g(x,u)+\a J\big(f(x,u)\big),$$
the control $u$ has the form $u=(u_1,\ldots,u_m)$, and each control component $u_\ell$ is constrained to take values in a finite set of real numbers $U_\ell(x)$. We assume that the cost per stage $g$ has the following strict monotonicity property for each $x\in X$ and $u_\ell\in U_\ell(x)$, $\ell=1,\ldots,m$:
$$g(x,u_1,\ldots,u_m)< g(x,u_1,\ldots,u_{\ell-1},\xi,u_{\ell+1},\ldots,u_m),\qquad \hbox{for all }\xi>u_\ell,\ \xi\in U_\ell(x).$$
We also assume that the system equation $x_{k+1}=f(x_k,u_k)$ is linear of the form
$$x_{k+1}=A x_k+Bu_k,$$
where $A$ and $B$ are matrices of appropriate dimension and have nonnegative components. The assumptions guarantee that for every $\m$, $J_\m$ is monotonically nondecreasing in $x$, and this this sufficient to that the function
$$H(x,u,J_\m)=g(x,u)+\a J_\m(Ax+Bu),$$
is strictly monotonically increasing along each control component $u_\ell$.
}
 
\vskip-1pc

\section{Concluding Remarks}
\vskip-0.5pc
\pn We have shown that in the  context of multiagent problems, agent-by-agent versions of the VI algorithm and related optimistic PI algorithms have greatly reduced computational requirements, while still maintaining a meaningful convergence property. While these algorithms may terminate with a suboptimal policy that is agent-by-agent optimal, they can be dramatically more efficient than the standard VI and optimistic PI algorithms, which may be computationally intractable even for a moderate number of agents. 

Several unresolved questions remain regarding algorithmic variations and conditions that guarantee that our algorithms obtain an optimal policy rather than one that is agent-by-agent optimal. Approximate versions of our algorithms of the type used in neuro-dynamic programming/reinforcement learning are also of interest, and are a subject for further investigation. Moreover, the basic idea of our approach, simplifying the minimization defining the VI operator while maintaining some form of convergence guarantee, can be extended in other directions to exploit special problem structures.

\vskip-1pc

\section{References}
\vskip-0.9pc
\def\ref{\vskip1.pt\pn}

\ref[BaP08] Bauso, D., and Pesenti, R., 2008.\ ``Generalized Person-by-Person Optimization in Team Problems with Binary Decisions," Proc.\ 2008 American Control Conference, pp.\ 717-722.

\ref[BaP12] Bauso, D., and Pesenti, R., 2012.\ ``Team Theory and Person-by-Person Optimization with Binary Decisions," SIAM Journal on Control and Optimization, Vol.\ 50, pp.\ 3011-3028.

\ref[Ber73] Bertsekas, D.\ P., 1973.\ ``Linear Convex Stochastic Control Problems Over an Infinite Horizon," IEEE Transactions on Aut.\ Control, Vol.\ AC-18, pp.\ 314-315.

\ref [BeT89]  Bertsekas, D.\ P., and Tsitsiklis, J.\ N., 1989.\ Parallel and Distributed Computation: Numerical Methods, Prentice-Hall, Englewood Cliffs, NJ; republished in 1996 by Athena Scientific, Belmont, MA.

\ref [BeT96]  Bertsekas, D.\ P., and Tsitsiklis, J.\ N., 1996.\ Neuro-Dynamic Programming, Athena Scientific, Belmont, MA.

\ref[BeY10] Bertsekas, D.\ P., and Yu, H., 2010.\ ``Asynchronous Distributed Policy Iteration in Dynamic Programming,"  Proc.\ of Allerton Conf.\ on Communication, Control and Computing,  Allerton Park, Ill, pp.\ 1368-1374.

\ref[BeY12] Bertsekas, D.\ P., and Yu, H., 2012.\ ``Q-Learning and Enhanced Policy Iteration in Discounted 
Dynamic Programming,"  Math.\ of OR, Vol.\ 37, pp.\ 66-94.

\ref[Ber10] Bertsekas, D.\ P., 2010.\
``Williams-Baird Counterexample for Q-Factor Asynchronous Policy Iteration,"
http://web.mit.edu/dimitrib/www/Williams-Baird Counterexample.pdf.

\ref[Ber12] Bertsekas, D.\ P., 2012.\ Dynamic Programming and Optimal Control, Vol.\ II, 4th edition, Athena Scientific, Belmont, MA.

\ref[Ber16] Bertsekas, D.\ P., 2016.\ Nonlinear Programming, 3rd Edition, Athena Scientific, Belmont, MA.

\ref[Ber18] Bertsekas, D.\ P., 2018.\ Abstract Dynamic Programming, Athena Scientific, Belmont, MA; on-line at http://web.mit.edu/dimitrib/www/RLbook.html.

\ref [Ber19]  Bertsekas, D.\ P., 2019.\ Reinforcement Learning and Optimal Control, Athena Scientific, Belmont, MA.

\ref[Ber20] Bertsekas, D.\ P., 2020.\ ``Multiagent Rollout Algorithms and Reinforcement Learning," arXiv preprint, arXiv:2002.07407.

\ref[Gup20] Gupta, A., 2020.\ ``Existence of Team-Optimal Solutions in Static Teams with Common Information: A Topology of Information Approach," SIAM J.\ on Control and Optimization, Vol.\ 58, pp.\ 998-1021.

\ref[Ho80] Ho, Y.\ C., 1980.\ ``Team Decision Theory and Information Structures," Proceedings of the IEEE, Vol.\ 68, pp.\ 644-654.

\ref[KLM82a] Krainak, J.\ L.\ S.\ J.\ C., Speyer, J., and Marcus, S., 1982.\ ``Static Team Problems - Part I: Sufficient Conditions and the Exponential Cost Criterion," IEEE Transactions on Automatic Control, Vol.\ 27, pp.\ 839-848.

\ref[KLM82b] Krainak, J.\ L.\ S.\ J.\ C., Speyer, J., and Marcus, S., 1982.\ ``Static Team Problems - Part II: Affine Control Laws, Projections, Algorithms, and the LEGT Problem," IEEE Transactions on Automatic Control, Vol.\ 27, pp.\ 848-859.

\ref[LTZ19] Li, Y., Tang, Y., Zhang, R., and Li, N., 2019.\ ``Distributed Reinforcement Learning for Decentralized Linear Quadratic Control: A Derivative-Free Policy Optimization Approach," arXiv preprint arXiv:1912.09135.

\ref[Mar55] Marschak, J., 1975.\ ``Elements for a Theory of Teams," Management Science, Vol.\ 1, pp.\ 127-137.

\ref[NMT13] Nayyar, A., Mahajan, A. and Teneketzis, D., 2013.\ ``Decentralized Stochastic Control with Partial History Sharing: A Common Information Approach," IEEE Transactions on Automatic Control, Vol.\ 58, pp.\ 1644-1658.

\ref[NaT19] Nayyar, A. and Teneketzis, D., 2019.\ ``Common Knowledge and Sequential Team Problems," IEEE Transactions on Automatic Control, Vol.\ 64, pp.\ 5108-5115.

\ref[Rad62] Radner, R., 1962.\ ``Team Decision Problems," Ann.\ Math.\ Statist., Vol.\ 33, pp.\ 857-881.

\ref[WiB93] Williams, R.\ J., and Baird, L.\ C., 1993.\ ``Analysis of
Some Incremental Variants of Policy Iteration: First Steps Toward
Understanding Actor-Critic Learning Systems,'' Report NU-CCS-93-11,
College of Computer Science, Northeastern Univ.,
Boston, MA. 

\ref[WaS00] de Waal, P.\ R., and van Schuppen, J.\ H., 2000.\ ``A Class of Team Problems with Discrete Action Spaces: Optimality Conditions Based on Multimodularity," SIAM J.\ on Control and Optimization, Vol.\ 38, pp.\ 875-892.

\ref[Wit68] Witsenhausen, H., 1968.\ ``A Counterexample in Stochastic Optimum Control," SIAM Journal on Control, Vol.\ 6, pp. 131?147.

\ref[Wit71a] Witsenhausen, H., 1971.\ ``On Information Structures, Feedback and Causality," SIAM J.\ Control, Vol.\ 9, pp.\ 149-160.

\ref[Wit71b] Witsenhausen, H., 1971.\ ``Separation of Estimation and Control for Discrete Time Systems," Proceedings of the IEEE, Vol.\ 59, pp.\ 1557-1566.

\ref[Wit88] Witsenhausen, H., 1988.\ ``Equivalent Stochastic Control Problems," Math.\ Control Signals Systems, Vol.\ 1, pp.\ 3-11.

\old{
\ref [WuB01] Wu, C., and Bertsekas, D.\ P., 2001.\ ``Distributed Power Control Algorithms for 
Wireless Networks,"  IEEE Trans.\ on Vehicular Technology, Vol.\ 50, pp.\ 504-514.
}

\ref[YuB13] Yu, H., and Bertsekas, D.\ P., 2013.\ ``Q-Learning and Policy Iteration Algorithms for Stochastic Shortest Path Problems," Annals of Operations Research, Vol.\ 208, pp.\ 95-132.

\ref[ZPB19] Zoppoli, R., Parisini, T., Baglietto, M., and Sanguineti, M., 2019.\ Neural Approximations for Optimal Control and Decision, Springer.

\end